\documentclass[12pt,leqno]{article}

\newtheorem{theorem}{Theorem}
\newtheorem{lemma}[theorem]{Lemma}
\newtheorem{proposition}[theorem]{Proposition}
\newtheorem{sublemma}[theorem]{Sublemma}
\newtheorem{definition}[theorem]{Definition}
\newtheorem{corollary}[theorem]{Corollary}
\newtheorem{problem}[theorem]{Problem}
\newtheorem{remark}[theorem]{Remark}
\newtheorem{claim}[theorem]{Claim}
\newtheorem{assumptions}[theorem]{Assumptions}
\newtheorem{examples}[theorem]{Examples}
\newtheorem{question}[theorem]{Question}
\newtheorem{sassumptions}[theorem]{Standing Assumptions}
\newtheorem{sassumption}[theorem]{Standing Assumption}
\newtheorem{conjecture}[theorem]{Conjecture}

\newcommand{\begintheorem}{\addtocounter{equation}{1}\begin{theorem}}
\newcommand{\beginlemma}{\addtocounter{equation}{1}\begin{lemma}}
\newcommand{\beginproposition}{\addtocounter{equation}{1}\begin{proposition}}
\newcommand{\beginsublemma}{\addtocounter{equation}{1}\begin{sublemma}}
\newcommand{\begindefinition}{\addtocounter{equation}{1}\begin{definition}}
\newcommand{\begincorollary}{\addtocounter{equation}{1}\begin{corollary}}
\newcommand{\beginproblem}{\addtocounter{equation}{1}\begin{problem}}
\newcommand{\beginremark}{\addtocounter{equation}{1}\begin{remark}}
\newcommand{\beginclaim}{\addtocounter{equation}{1}\begin{claim}}
\newcommand{\beginassumptions}{\addtocounter{equation}{1}\begin{assumptions}}
\newcommand{\beginexamples}{\addtocounter{equation}{1}\begin{examples}}
\newcommand{\beginquestion}{\addtocounter{equation}{1}\begin{question}}
\newcommand{\beginsassumptions}{\addtocounter{equation}{1}\begin{sassumptions}}
\newcommand{\beginsassumption}{\addtocounter{equation}{1}\begin{sassumption}}
\newcommand{\beginconjecture}{\addtocounter{equation}{1}\begin{conjecture}}

\begin{document}

\title{Some topics concerning analysis on metric spaces and semigroups
of operators}

\author{Stephen Semmes}

\date{}

\maketitle

\subsubsection*{Classical analysis on Euclidean spaces (as in \cite{St2, SW})}

	Fix a positive integer $n$, and consider the Euclidean space
${\bf R}^n$ equipped with the standard distance function $|x-y|$
and Lebesgue measure.  If $f(x)$ is a locally-integrable function
on ${\bf R}^n$, then the \emph{Hardy--Littlewood maximal function}
$f^*(x)$ associated to $f$ is defined by
\begin{equation}
	f^*(x) = \sup_{x \in B} \frac{1}{\mu(B)} \int_B |f(y)| \, d\mu(y),
\end{equation}
where the supremum is taken over all open balls $B$ in ${\bf R}^n$
which contain $x$.  The supremum may be $+ \infty$, so that $f^*$ is
actually a (nonnegative) extended real-valued function.  This is
sometimes referred to as the \emph{uncentered} maximal function, and
there are variants defined in terms of balls centered at $x$, or using
cubes instead of balls.  A nice feature of $f^*(x)$ is that it is
upper semicontinuous, which is to say that
\begin{equation}
	\{x \in {\bf R}^n : f^*(x) > t\}
\end{equation}
is an open subset of ${\bf R}^n$ for each positive real number $t$.
Indeed, if $f^*(x) > t$ for some $x$, then there is a ball $B$ containing
$x$ such that 
\begin{equation}
	\frac{1}{\mu(B)} \int_B |f(y)| \, d\mu(y) > t,
\end{equation}
and it follows that $f^*(z) > t$ for all $z$ in $B$.

	Clearly the supremum of $f^*$ is less than or equal to the
$L^\infty$ norm of $f$.  A famous \emph{weak type $(1,1)$} result says
that 
\begin{equation}
	|\{x \in M : f^*(x) > t\}| \le C(n) \, t^{-1} \|f\|_1,
\end{equation}
for some constant $C(n) > 0$ that depends only on the dimension $n$
and all functions $f$, where $|E|$ denotes the Lebesgue measure of a
set $E$ and $\|f\|_1$ is the usual $L^1$ norm of $f$.  In particular,
$f^*$ is finite almost everywhere in this case.  For $p > 1$ there
is a \emph{strong type} result, which means that
\begin{equation}
	\|f^*\|_p \le C(n,p) \, \|f\|_p
\end{equation}
for some constant $C(n,p)$ which depends only on $n$ and $p$, and
where $\|f\|_p$ denotes the usual $L^p$ norm of $f$.  This can in fact
be derived from the preceding estimates for $p = 1, \infty$ through
a general interpolation result.

	One might be interested in other kinds of averages of $f$,
such as those given by integrating $f$ against the Poisson kernel
or the Gauss--Weierstrass kernel.  These are exactly the quantities
which arise in the extensions of $f$ to the upper half space
${\bf R}^n \times (0, \infty)$ which are harmonic or satisfy the
heat equation (and which satisfy additional mild growth conditions
to avoid modest ambiguities).  Fortunately, these averages
can be estimated in terms of averages over balls in a simple
way, so that the corresponding maximal functions are bounded in
terms of $f^*$.  Thus the inequalities above for $f^*$ provide basic
results about the boundary behavior of solutions to the Laplace and
heat equations on ${\bf R}^n \times (0, \infty)$.

	Some other interesting operators are the singular integral
operators
\begin{equation}
	R_j(f)(x) = p.v. \int_{{\bf R}^n} \frac{x_j - y_j}{|x - y|^{n+1}}
					\, f(y) \, dy,
\end{equation}
$1 \le j \le n$, and
\begin{equation}
	I_{it}(f)(x) = p.v. \int_{{\bf R}^n} \frac{1}{|x - y|^{n + i t}}
					\, f(y) \, dy,
\end{equation}
$t \in {\bf R}$, $t \ne 0$.  Some care is involved in taking the
principal values, especially in the second case.  For $I_{it}$,
different ways of defining the principal values will even lead to
different answers, but the difference is rather mild (a multiple of
the identity operator).

	One can show that these operators are bounded on $L^2$ using
special structure related to $p = 2$, i.e., Fourier transform and
Hilbert space methods.  This can be extended to boundedness on $L^p$
when $1 < p < \infty$ and the weak type $(1,1)$ property for $p = 1$
using well-known techniques in harmonic analysis.  For $p = \infty$
there are estimates in terms of BMO, as a substitute for $L^\infty$
bounds which do not work.  Similar results apply to numerous other
operators of similar type.

\subsubsection*{Spaces of homogeneous type \cite{CW1, CW2}}

	A \emph{space of homogeneous type} can be described as a
triple $(M, d(x,y), \mu)$, where $M$ is a nonempty set, $d(x,y)$ is a
metric on $M$ (and thus is a symmetric nonnegative real-valued
function on $M \times M$ which vanishes exactly when $x = y$ and
satisfies the triangle inequality), and $\mu$ is a \emph{doubling
measure} on $M$.  The latter means that $\mu$ is a nonnegative Borel
measure which assigns positive finite measure to open balls in $M$,
and for which there is a constant $C > 0$ such that
\begin{equation}
\label{doubling condition for mu}
	\mu(B(x,2r)) \le C \, \mu(B(x,r))
\end{equation}
for every ball $B(x,r)$ in $M$.  Of course (\ref{doubling condition
for mu}) implies that $\mu$ assigns positive finite measure to every
open ball in $M$ as soon as this holds for a single such ball.  One
might also ask that $M$ be complete, in the sense that Cauchy
sequences converge, or that open subsets of $M$ be realizable as
countable unions of compact sets.  Basic examples of spaces of
homogeneous type are given by Euclidean spaces with the standard
metric and Lebesgue measure.  Reasonably-smooth domains or manifolds
are also included in this notion.

	It can be convenient to allow $d(x,y)$ to be a
\emph{quasimetric} instead of a metric, which means that a positive
constant factor is allowed on the right side of the triangle
inequality, and the notion of a space of homogeneous type is often
formulated in this manner.  As in \cite{MS1}, there are always metrics
not too far from quasimetrics, so that for many purposes one might as
well restrict to metrics.

	The Hardy--Littlewood maximal function $f^*$ associated to a
locally-integrable function $f$ can be defined on a space of
homogeneous type in the same manner as on Euclidean spaces.  A basic
fact is that the weak type $(1,1)$ estimate extends to this general
setting.  The supremum of $f^*$ is still bounded by the $L^\infty$
norm of $f$, and $L^p$ estimates for $1 < p < \infty$ follow from the
$p = 1, \infty$ estimates through general interpolation arguments, as
before.

	Another basic result is that one has ``Calder\'on--Zygmund
inequalities'' for singular integral operators analogous to those on
${\bf R}^n$.  That is, one can start with a linear operator $T$ which
is bounded on $L^2$, or some other fixed $L^{p_1}$, and which is
associated to a kernel that satisfies suitable size and smoothness
conditions, and derive boundedness on $L^p$ for all $1 < p < \infty$
and a weak-type inequality for $p = 1$.  One can also get BMO
extimates for $p = \infty$, estimates on Hardy spaces as an
alternative to the weak-type inequality for $p = 1$ as well as
allowing for some $p < 1$, etc.  The compatibility between the metric
and the measure given by the doubling condition is quite remarkable.

	Let us mention two classes of examples of spaces of
homogeneous type which were examined on their own before the general
notion.  In the first case, which was studied by my colleague Frank
Jones \cite{Frank}, one takes ${\bf R}^n \times {\bf R}$ with the
distance between two points $(x,s)$, $(y,t)$ defined to be
\begin{equation}
	|x-y| + |s-t|^{1/2},
\end{equation}
where $|x-y|$, $|s-t|$ denote the usual distances in ${\bf R}^n$,
${\bf R}$, respectively.  Sometimes other expressions are used for
essentially the same geometry; a key point is that the distance
behaves well under the non-isotropic dilations 
\begin{equation}
	(x,t) \mapsto (r x, r^2 t),
\end{equation}
for $r > 0$, just as the ordinary metric on ${\bf R}^n$ behaves well
under the dilations $x \mapsto r x$.  Of course the metric is also
invariant under translations on ${\bf R}^n \times {\bf R}$, and is
compatible with the usual topology.  For the measure one still uses
Lebesgue measure.  The measure of a ball of radius $\rho$ is a
constant multiple of $\rho^{n+2}$, and the doubling condition is
satisfied.  In this case the singular integral theory can be applied
to operators related to the heat operator, whereas the standard
geometry on ${\bf R}^n$ fits with operators related to the Laplacian.
Note that there is a kind of tricky point here, in which the $t$ parameter
is included in the underlying space.  

	A second basic situation corresponds to the unit sphere in
${\bf C}^n$, which, for $n \ge 2$, has a non-Euclidean geometry which
is adapted to several complex variables, holomorphic functions on the
unit ball in ${\bf C}^n$, etc.  Just as in the previous case, one can
still use ordinary Lebesgue measure on the sphere, and this measure is
doubling with respect to the non-Euclidean geometry.  (For that
matter, it is also doubling with respect to the usual Euclidean
geometry.)  The Hardy--Littlewood maximal function with respect to the
non-Euclidean geometry is closely connected to maximal functions and
limits for holomorphic functions in the ball along certain
``admissible'' regions, just as the classical maximal function is
connected to nontangential maximal functions for holomorphic functions
in one complex variable or harmonic functions in several real
variables.  A fundamental singular integral operator in this situation
is the \emph{Szeg\"o projection}, which is the orthogonal projection
from $L^2$ of the unit sphere onto the subspace of functions which are
boundary values of holomorphic functions on the ball.  This operator
is bounded on $L^2$ with norm $1$ by definition, and its kernel can be
computed explicitly.  With respect to the non-Euclidean geometry, the
kernel satisfies the appropriate size and smoothness conditions, so
that the operator is in fact bounded on $L^p$, $1 < p < \infty$, and
so on.  See \cite{Kor1, Kor2, Kor3, KV1, KV2, Krantz2, Rudin, St3,
St4}.

\subsubsection*{Semigroups of operators}

	In another direction, suppose that $\mathcal{B}$ is a Banach
space, and that $\{T_t\}_{t \ge 0}$ is a \emph{semigroup of bounded
operators} on $\mathcal{B}$.  Specifically, assume that $T_0$ is the
identity operator $I$, that the operator norm of $T_t$ is bounded by
some constant $k$ for $0 \le t \le 1$, that
\begin{equation}
	T_{s + t} = T_s \circ T_t
\end{equation}
for all $s, t \ge 0$, and that $\lim_{t \to 0} T_t(f) = f$ for all $f$
in $\mathcal{B}$.  Of course the semigroup property together with the
uniform bound for the operator norm of the $T_t$'s for $0 \le t \le 1$
implies an exponentially-increasing bound for the operator norm of
$T_t$ for all $t$'s.

	There is a remarkable amount of mathematics around this kind
of situation.  In fact, this is just the beginning; one can add
relatively-simple hypotheses which occur in numerous settings and
which add quite a bit more structure.  As a basic distinction, one
might think of $\{T_t\}_{t \ge 0}$ as being a semigroup of unitary
transformations on a Hilbert space, or a semigroup of invertible
linear mappings on a Banach space more generally, as is associated to
solutions of a wave equation, or one might think of $\{T_t\}_{t \ge
0}$ as defining a diffusion, as is associated to solutions of a heat
equation.

	Here we shall mostly focus on the second type of situation.
We assume now that we have a measure space $M$ with a positive measure
$\mu$, and we take for our Banach space $\mathcal{B}$ the Hilbert
space $L^2(M, \mu)$.  We ask too that each $T_t$ be self-adjoint and
\emph{positivity-preserving}, which means that for each nonnegative
function $f$ on $M$, $T_t(t)$ is also a nonnegative function on $M$
for every $t \ge 0$.  Each of these conditions is significant in its
own right, and part of the beauty of the subject arises from the
interplay between them.

	Let us also ask that the $T_t$'s extend to bounded operators
on $L^p(M, \mu)$ for each $1 \le p \le \infty$, and in fact that the
$T_t$'s are \emph{contractions} on all $L^p$, which is to say that the
operator norms are all less than or equal to $1$.  If ${\bf 1}$
denotes the function on $M$ which is identically equal to $1$, then we
ask that $T_t({\bf 1}) = {\bf 1}$ for all $t$.  These conditions
are satisfied by the semigroups associated to the heat kernel and
Poisson kernel on ${\bf R}^n$, for instance.

	A famous result of Stein states that the maximal function
inequalities
\begin{equation}
	\|\sup_{t > 0} |T_t(f)| \|_p \le A_p \, \|f\|_p
\end{equation}
hold for $1 < p \le \infty$, i.e., for some constant $A_p$ and all
functions $f$ in $L^p$.  In this setting there is also a ``singular
integral operator'' theory, for operators which are functions of the
generator of the semigroup.  Boundedness on $L^2$ for these operators
is easily determined through the spectral representation.  Some
general conditions for boundedness on $L^p$ are described in
\cite{St1}, and of course more precise information depends on the
particular situation.

	The significance of self-adjointness is illustrated by the
example where $T_t$ is defined on functions on ${\bf R}$ to be
translation by $t$.  This semigroup of operators is
positivity-preserving and preserves all $L^p$ norms, but the maximal
inequality fails completely for $p < \infty$.  In this case it is
natural to consider averages of $T_t f$ and suprema of the averages,
as in ergodic theory, and as for semigroups associated to
measure-preserving transformations on the underlying measure space
more generally.

\subsubsection*{Semigroups and geometry}

	There are very interesting combinations of the spaces of
homogeneous type and semigroups of operators pictures, involving
bounds for kernels of semigroups, and $L^p$ mapping properties of
operators related to the semigroup.  See \cite{AuMT, AuT, CD1, CD2,
CD3, DM, DR}, for instance.  Another perspective has recently been
studied in \cite{GHL}, with the following set-up.  One assumes again
that $(M, d(x,y))$ is a metric space, that $\mu$ is a positive Borel
measure on it, and that $T_t$ is a symmetric contraction semigroup of
linear operators as before.  Now one asks in addition that for $t > 0$
the operator $T_t$ is defined by a nonnegative kernel $k_t(x,y)$, so
that
\begin{equation}
	T_t(f)(x) = \int_M k_t(x,y) \, f(y) \, d\mu(y).
\end{equation}
For the kernel $k_t(x,y)$ one considers upper and lower bounds of the
form
\begin{equation}
   \frac{1}{t^{\alpha/\beta}} \, \phi_1\biggl(\frac{d(x,y)^\beta}{t}\biggr)
		\le k_t(x,y) \le
   \frac{1}{t^{\alpha/\beta}} \, \phi_2\biggl(\frac{d(x,y)^\beta}{t}\biggr).
\end{equation}
Here $\alpha$, $\beta$ are positive constants, and $\phi_1(u)$,
$\phi_2(u)$ are monotone decreasing positive functions on
$[0,\infty)$, with $\phi_1(u_1) > 0$ for some $u_1 > 0$ and
$\phi_2(u)$ normally asked to satisfy decay conditions.

	The parameter $\alpha$ is related to volume growth in $M$, and
this is discussed in \cite{GHL}.  The connection between $\beta$ and
the geometry of $M$ is also treated in \cite{GHL}.  For the standard
heat semigroup on Euclidean spaces, $\beta$ is always equal to $2$.
There are heat semigroups associated to subelliptic operators in place
of the ordinary Laplacian which also satisfy these conditions with
$\beta = 2$, with respect to an associated metric.  A basic version of
this arises for the unit sphere in ${\bf C}^n$, $n \ge 2$, and
non-Euclidean geometry on it, as indicated earlier.  There are a
number of fractals such as Sierpinski gaskets and carpets and
semigroups on them which satisfy the conditions above with various
values of $\beta$.  Compare with \cite{BB1, BB2, BB3, BB4, BP, FHK,
Kigami, Kumagai}.

	Without the semigroup property, there are well-known fairly
simple constructions of approximations to the identity on spaces of
homogeneous type with nice properties.  The semigroup property of
course imposes very strong restrictions.  For that matter,
commutativity of the operators in the family is a substantial
condition.

	Let us note that decay conditions on $\phi_2(u)$ above can be
quite significant.  A very nice contraction semigroup on ${\bf R}^n$
is given by the Poisson kernel, and for this kernel the decay is not
very fast.  Modest decay conditions for the kernel are adequate for a
number of applications, even if they are not sufficient for other
results, as in \cite{GHL}.

\subsubsection*{Analysis on fractals like Sierpinski gaskets and carpets}

	Of course decay conditions for the kernel of a semigroup are
closely connected to locality conditions for the generator of the
semigroup.  In the classical cases on Euclidean spaces or tori for
periodic functions, etc., the heat kernel has fast decay and is
generated by the Laplace operator, while the Poisson kernel does not
have very fast decay and is generated by a constant multiple of the
\emph{square root} of the Laplace operator, which is not a local
operator.

	The fast decay for the kernels of the semigroups on Sierpinski
gaskets and carpets mentioned before reflects the fact that the
generators are nice local operators, versions of ``Laplacians'' for
these fractals.  The fractal structures play a role and are reflected
in the parameter $\beta$, but still there are nice operators which are
like differential operators.

	For example, there are remarkable results concerning elliptic
and parabolic Harnack inequalities for these operators.  See
\cite{BB1, BB2, BB3, BB4, Kigami}.

\end{document}